\documentclass{amsart}
\usepackage{amscd}
\numberwithin{equation}{section}
\usepackage{amssymb}
\let\cal\mathcal
\def\Ascr{{\cal A}}

\def\Cscr{{\cal C}}
\def\Dscr{{\cal D}}
\def\Escr{{\cal E}}
\def\Fscr{{\cal F}}

\def\Lscr{{\cal L}}
\def\Mscr{{\cal M}}

\def\Oscr{{\cal O}}

\def\Qscr{{\cal Q}}

\def\Uscr{{\cal U}}

\let\blb\mathbb
 
\def\DD{{\blb D}}

\def \PP{{\blb P}}

\def \ZZ{{\blb Z}}

\def \NN{{\blb N}}

\def\udim{\operatorname{\underline{\dim}}}
\def\Sym{\operatorname{Sym}}

\def\Id{\operatorname{id}}
\def\pr{\mathop{\text{pr}}\nolimits}

\def\Lboxtimes{\overset{L}{\boxtimes}}
\def\Lotimes{\overset{L}{\otimes}}
\def\quot{/\!\!/}

\def\Mod{\operatorname{Mod}}
\def\mod{\operatorname{mod}}
\def\Gr{\operatorname{Gr}}

\def\gr{\operatorname{gr}}

\def\Qch{\operatorname{Qch}}
\def\coh{\mathop{\text{\upshape{coh}}}}

\def\rad{\operatorname {rad}}

\def\gr{\operatorname {gr}}
\def\Spec{\operatorname {Spec}}

\def\GL{\operatorname {GL}}

\def\diag{\operatorname {diag}}
\def\Ext{\operatorname {Ext}}
\def\Hom{\operatorname {Hom}}

\def\End{\operatorname {End}}
\def\RHom{\operatorname {RHom}}
\def\uRHom{\operatorname {R\mathcal{H}\mathit{om}}}
\def\Sl{\operatorname {Sl}}

\def\ker{\operatorname {ker}}

\def\End{\operatorname {End}}

\def\rk{\operatorname {rk}}

\def\gldim{\operatorname {gl\,dim}}

\def\r{\rightarrow}

\DeclareMathOperator{\Proj}{Proj}

\DeclareMathOperator{\Tors}{Tors}

\DeclareMathOperator{\HEnd}{\mathcal{E}\mathit{nd}}

\let\dirlim\injlim

\newtheorem{lemma}{Lemma}[section]
\newtheorem{proposition}[lemma]{Proposition}
\newtheorem{theorem}[lemma]{Theorem}

\newtheorem{lemmas}{Lemma}[subsection]
\newtheorem{propositions}[lemmas]{Proposition}
\newtheorem{theorems}[lemmas]{Theorem}

\theoremstyle{definition}

\newtheorem{example}[lemma]{Example}
\newtheorem{definition}[lemma]{Definition}
\newtheorem{conjecture}[lemma]{Conjecture}
\newtheorem{question}[lemma]{Question}
\newtheorem{examples}[lemmas]{Example}

{

}

\theoremstyle{remark}

\newtheorem{remark}[lemma]{Remark}
\newtheorem{remarks}[lemmas]{Remark}

\newdimen\uboxsep \uboxsep=1ex
\def\uboxn#1{\vtop to 0pt{\hrule height 0pt depth 0pt\vskip\uboxsep
\hbox to 0pt{\hss #1\hss}\vss}}

\def\uboxs#1{\vbox to 0pt{\vss\hbox to 0pt{\hss #1\hss}
\vskip\uboxsep\hrule height 0pt depth 0pt}}

\author{Michel Van den Bergh}
 \address{Universiteit Hasselt\\ Universitaire Campus\\ 3590 Diepenbeek}
  \email{michel.vandenbergh@uhasselt.be}
\thanks{The author is a director of research at the FWO}
\date{Oct 25, 2002}
\keywords{Resolutions, derived equivalence, homologically homogeneous rings}
\subjclass{Primary 18E30, 14E30, 14A22}
\title{Non-commutative crepant resolutions\\(with some corrections)}

\begin{document}
\begin{abstract}
We introduce the notion of a ``non-commutative crepant'' resolution
of a singularity and show that it exists in certain cases. We also
give some evidence for an extension of a conjecture by
Bondal and Orlov, stating that different crepant resolutions
of a Gorenstein singularity have the same derived category.
\end{abstract}
\maketitle
\tableofcontents
\section{Introduction}
Let $k$ be an algebraically closed field and let $R$ be an integral
Gorenstein $k$-algebra of dimension $n$.
Put $X=\Spec R$ and let $f:Y\r X$ be a resolution of singularities.
$f$ is a crepant resolution of $X$ if $f^\ast\omega_X=\omega_Y$.
Crepant resolutions do not always exist and they are usually not
unique. A conjecture of Bondal and Orlov states that if $f_1:Y_1\r X$,
$f_2:Y_2\r X$ are different crepant resolutions of $X$ then
$D^b(\coh(Y_1))\cong D^b(\coh(Y_2))$.

In this paper we study certain non-commutative analogues of crepant
resolutions and we call these ``non-commutative crepant resolutions''.
A non-commutative crepant resolution of $R$ is an algebra
$A=\End_R(M)$ where $M$ is a reflexive $R$-module and where $A$ has
finite global dimension and is a (maximal) Cohen-Macaulay $R$-module. For 
some general discussion of non-commutative resolutions we refer to
\cite[\S5]{BO1}

 A
standard example of a non-commutative crepant resolution is the following:
\begin{example} Let $G$ be a finite group and let $V$ be a finite dimensional
  $G$-representation such that $G\subset \Sl(V)$. Put $S=\Sym(V)$ and
  $R=S^G$. Then ${{A}}=\End_R(S)\cong S\ast G$ is a non-commutative
  crepant resolution of $R$.
\end{example}
In \cite{VdB31} we made the following conjecture.
\begin{conjecture}
\label{ref-1.2-0}
If $R$ is three dimensional and has terminal singularities
then it has a non-commutative crepant resolution if and only if it has
a crepant commutative resolution.
\end{conjecture}
We prove this conjecture below (Theorem \ref{ref-6.6.3-38}). We also
show that a version of the Bondal-Orlov conjecture holds in the sense
that the two resolutions in Conjecture \ref{ref-1.2-0} are derived
equivalent.
One direction in the proof of Conjecture \ref{ref-1.2-0}
is obtained from the beautiful version of the McKay correspondence
 given in \cite{BKR} (with virtually the same proof).

In addition to proving Conjecture \ref{ref-1.2-0} we will give
two other instances where a non-commutative resolution exists:
\begin{enumerate}
\item Cones over Del Pezzo surfaces (see~\S\ref{ref-7-39}).
\item Invariants of a one  dimensional torus acting linearly on a polynomial
ring (see~\S\ref{ref-8-46}).
\end{enumerate}
We would like to thank Paul Smith and Amnon Yekutieli for some useful
discussions. 

\medskip

The main reason for this new version is that the argument for the
existence of non-commutative crepant resolutions for cones of Del
Pezzo surfaces was incorrect in the published version of this paper
\cite{VdB32}. Luckily the statement follows easily from
the work of Kuleshov and Orlov.  This approach was suggested to the
author by Tom Bridgeland.

Besides this we have updated a few references but we have made no
attempt to cover recent developments (this would have required
extensive changes to the paper).

Finally we have corrected a few other minor errors. 

\begin{itemize}
\item The proof of Proposition \ref{ref-3.3-5} was slightly faulty.
\item There was a disastrous typo in the statement of Lemma
  \ref{wellknown}(5)  (pointed out by Orlov).
\end{itemize}
\section{Conventions}
\label{ref-2-1}
Throughout $k$ is an algebraically closed field of characteristic
zero. All rings are $k$-algebras. Modules are right modules and
Cohen-Macaulay means maximal Cohen-Macaulay.  When we say that $R$ is
complete local then we mean that $R$ is a (commutative) complete local
$k$-algebra with residue field $k$.

Below when we say ``reduction to the complete case'' we usually mean
completing the strict Henselization of a commutative noetherian ring
at a prime ideal and then lifting the residue field. Replacing $k$
with this lifted field we are in a ``complete local'' situation,
in the sense of the previous paragraph.

We also use it in the following sense: let $R$ be a commutative
noetherian ring and let $R\r K$ be a map of $R$ to a field with kernel
$P$. Then a complete local situation may be obtained by replacing $R$
with a suitable flat complete local extension of the completion of $R$ at $P$
which has residue field $K$ and which also contains $K$.

Below we will sometimes use the theory of dualizing complexes for schemes
and for non-commutative coherent sheaves of rings over these.  When
the underlying scheme is of finite type or when it is the spectrum of a
complete local ring then there is a canonical dualizing complex which
we will call the ``Grothendieck dualizing complex'' and which we
denote by $D$.  One possible way of characterizing this dualizing
complex is through the theory of ``rigid dualizing complexes''
\cite{VdB16,YZ1,YZ3} (and its topological variant). The duality
functor associated to $D$ will be denoted by $\mathbb{D}$.

We frequently use the following result by Keller  \cite[Thm 4.3]{Keller1}:
\begin{theorem}
\label{ref-2.0.1-2}
Let $\Escr$ be a Grothendieck category and
  assume that $\Ascr=D(\Escr)$ is generated by a compact
  object $E$ (i.e. $E^\perp=0$ and $\Hom_\Ascr(E,-)$ commutes with direct
sums). Then $\Ascr=D(\Lambda)$ where $\Lambda$ is a
  DG-algebra whose cohomology is given by $\Ext^\ast(E,E)$.
\end{theorem}
In our case we will always have $\Ext^i(E,E)=0$ for $i\neq 0$ and hence
$\Lambda$ is a true algebra. In that case we refer to $E$ as a \emph{tilting
object}.
\section{Background on homologically homogeneous rings}
\label{ref-3-3}
It is well-known that for a non-commutative ring the condition of
being of finite global dimension is very weak. In order to obtain
a good homological theory generalizing the commutative case, a stronger
condition is needed. Such a condition was introduced in \cite{BH}.

Let $R$ be a commutative noetherian $k$-algebra and let ${{A}}$
be a module-finite $R$-algebra. We will say that ${{A}}$ is 
\emph{homologically
homogeneous} if for all $P\in \Spec R$ we have
\begin{itemize}
\item[(H1)] $\gldim {{A}}_P=\dim R_P$.
\item[(H2)] ${{A}}_P$ is (maximal) Cohen-Macaulay.
\end{itemize}
If $R$ is equidimensional of dimension $n$ 
this is is equivalent to the following condition given in \cite{BH}.
\begin{itemize}
\item[(H1')] All simple ${{A}}$-modules have the same projective dimension
$n$.
\end{itemize}
\begin{example}
Consider the ring 
\[
{{A}}=\begin{pmatrix} R& m\\
R & R
\end{pmatrix}
\]
with $R=k[x,y]$, $m=(x,y)\subset R$. One has $\gldim {{A}}=2$ but ${{A}}$
is not Cohen-Macaulay.  So ${{A}}$ is \emph{not} homologically homogeneous.
\end{example}
Zero-dimensional h.h.\ rings are semi-simple. One-dimensional h.h.\ 
rings are classical hereditary orders \cite{reiner}. Two
dimensional h.h.\ rings over a complete local ring were classified in
\cite{Artinx,RVdB}. The results in \cite{AS,ATV1,ATV2,Steph1,Steph2} may be viewed as a
classification of three-dimensional graded local h.h.\ rings.

Not everything about homologically homogeneous rings is understood. For
example it
is an interesting question to understand precisely which rings can
occur as centers of homologically homogeneous rings. In particular we
would like to know the following:
\begin{question} 
\label{ref-3.2-4}  
If ${{A}}$ is  homologically homogeneous and $R$ is the
center of ${{A}}$,
is it  true that $R$ has rational singularities?
\end{question}
This question is answered affirmatively in complete generality in
\cite{VdBSt1}. Since we don't want to change too much we have kept the
discussion below although it is somewhat obsolete.  

\medskip

In the published version of this paper \cite{VdB32}
we were able to give a positive answer to Question \ref{ref-3.2-4} in
the graded case. Unfortunately as written the proof was not quite
correct. So below we give a corrected proof.
\begin{proposition} 
\label{ref-3.3-5}
Assume that ${{A}}$ is a finitely generated positively graded homologically
homogeneous 
graded $k$-algebra with semi-simple  degree zero part. Let $R$ be
the center of ${{A}}$ and assume that $R$ has an isolated
singularity. Then $R$ has rational singularities.
\end{proposition}
\begin{proof} 
The part of degree zero of $R$ is a direct sum of fields. Using the 
resulting central idempotents we may 
reduce to the case where $R$ is connected (i.e. $R_0=k$).
 Put $n=\dim R$.
If $n=0$ then there is nothing to prove so we assume $n>0$. 
Put  $m=R_{>0}$. $M={{A}}_{>0}$. Since $R$ has an isolated
singularity, it will have a rational singularity if and only if
the following condition holds  \cite{Watanabe}:
\begin{equation}
\label{ref-3.1-6}
\begin{split}
H^i_m(R)=0&\qquad\text{ for $i<n$}\\
H^n_m(R)_u=0&\qquad\text{ for $u\ge 0$}
\end{split}
\end{equation}

Since we are in characteristic zero, $R$ is a direct summand of ${{A}}$ 
and hence it is certainly Cohen-Macaulay. Thus $H^i_m(R)=0$ for $i<n$.

Let $S$ be a graded simple right ${{A}}$-module. 
From the theory of homologically homogeneous rings \cite{BH,VdBSt1} we
obtain $\Ext^n_{{A}}(S,{{A}})$ is a graded simple left $A$-module. Replacing $S$ by a minimal
(finite!) projective $A$-resolution we easily see that $\Ext^n_{{A}}(S,{{A}})$
is concentrated in strictly negative degree.

  From this we easily obtain
that $ H^n_{M}({{A}})=\dirlim \Ext^n_A({{A}}/M^n,{{A}}) $ lives in
strictly negative degree. Since $H^n_{M}({{A}})= H^n_m({{A}})$ (see
e.g. \cite{AZ}) and $H^n_m(R)$ is a direct summand of $H^n_m({{A}})$
we obtain that $H^n_m(R)$ also lives in strictly negative degree.
\end{proof}

\medskip

In \cite{RVdB} it is shown that every two-dimensional homologically homogeneous
$k$-algebra over a complete local ring is obtained as the completion
of a graded one. This is unfortunately false in higher dimension (otherwise
 Proposition \ref{ref-3.3-5} would yield an affirmative answer to
Question \ref{ref-3.2-4}).
\begin{example}
 We construct a three-dimensional homologically homogeneous $k$-algebra
which is not the completion of a graded one. Assume that $R$ is a complete
  local ring with a three-dimensional terminal Gorenstein singularity
  which is \emph{not} quasi-homogeneous 
 but which is such that $X=\Spec R$ has a crepant
  resolution. Then it follows from Theorem \ref{ref-5.1-12} below that there
  is a reflexive $R$-module such that ${{A}}=\End_R(M)$ is
  homologically homogeneous.  Clearly this ${{A}}$ is not the
  completion of a graded ring since otherwise the same would be true for its
center 
  $R$ also.

To make the example explicit consider
$
R=k[[x,y,z,t]]/(f)
$
where
\[
f=xy-(z-t^2)(z-t^3)(z-t^4)
\]
This is a compound Du Val singularity which has a crepant resolution. Since $R$
has a hypersurface singularity it will be quasi-homogeneous if and
only if the Milnor number $\mu(f)$ and the Tyurina number $\tau(f)$ of
$f$ are equal \cite{Saito}.
Using the computer algebra program SINGULAR \cite{SINGULAR} we compute:
$
\mu(f)=12$,
$\tau(f)=11$.
So
$R$ is \emph{not} quasi-homogeneous.
\end{example}
A more sensible approach to Question \ref{ref-3.2-4} is probably
through the following:
\begin{question} If ${{A}}$ is a homologically homogeneous ring 
over a complete local ring, is there  a separated filtration 
${{A}}={{A}}_0\supset {{A}}_1\supset\cdots$ 
 on ${{A}}$ 
with ${{A}}_1=\rad({{A}})$
(maybe
the $\rad({{A}})$-adic filtration(?)) such that $\gr{{A}}$ is still
homologically homogeneous? 
\end{question}
Since a deformation of a rational singularity is rational
\cite{Elkik}, it not hard to see that an affirmative answer to this
question combined with Proposition \ref{ref-3.3-5} would
yield an affirmative answer to Question \ref{ref-3.2-4}.
\begin{remark}
It is tempting to try to develop directly a theory of rational singularities
for non-commutative rings. For example \eqref{ref-3.1-6} seems to
make perfect sense for non-commutative graded rings. Unfortunately
it is easy to see that \eqref{ref-3.1-6} is not invariant under
Morita-equivalence so there are certainly some problems.

A related objection is that there seems to be no obvious
non-commutative analogue of the Grauert-Riemenschneider vanishing
theorem. It seems unlikely that there is a good theory of rational 
singularities without the Grauert-Riemenschneider theorem.
\end{remark}
\section{Non-commutative crepant resolutions}
\label{ref-4-7}
Unless otherwise specified, in the rest of this paper $R$ is a normal
Gorenstein domain and $X=\Spec R$. Let $f:Y\r X$ be a projective
morphism with $Y$ regular.  We
say that $f$ is \emph{crepant} if $f^\ast \omega_X=\omega_Y$. 

Bondal and Orlov
conjecture in \cite{Bondal1} that different crepant resolutions
of $X$ have equivalent bounded derived
categories of coherent sheaves. If $X$ is projective then it is known that
they have the same Hodge numbers \cite{Batyrev,Kont1}.

We propose a definition of a non-commutative crepant resolution. The 
motivation is mainly philosophical and the main example is given
by the McKay correspondence \cite{BKR} (see also \cite{KV} in dimension two).

\begin{definition} 
\label{ref-4.1-8}
\emph{A 
non-commutative
crepant} resolution of $R$ is an homologically homogeneous
$R$-algebra of the form ${{A}}=\End_R(M)$ 
where $M$ is a reflexive $R$-module.
\end{definition}
This definition will serve as a vehicle for introducing the various
examples we discuss in this paper. Some obvious variations are
possible (see Remark \ref{ref-4.3-9} below).

The following remark will be convenient.
\begin{lemma}
If $A=\End_R(M)$ for a reflexive $R$-module $M$ and if $\gldim A<\infty$ and
$A$ is Cohen-Macaulay then $A$ is homologically homogeneous (and hence
$A$ is a non-commutative crepant resolution).
\end{lemma}
\begin{proof} We may assume that $R$ is local. Put $n=\dim R$. 
Using the fact that $\Ext^i_R(A,R)=0$ for $i>0$  (since
$A$ is Cohen-Macaulay)  we easily deduce for $M\in \mod(A)$: $\Ext^i_A(M,A)=
\Ext^i_R(M,R)$. Thus the injective dimension of $A_A$ is $\le n$. Applying
the formula with $M$ a simple $A$-module obtain that the injective
dimension of $A_A$ is precisely $n$.

Since the global dimension of $A$ equal to the injective dimension of
$A_A$ (given that it is finite) we find $\gldim A\le n$.
\end{proof}

The point of Definition \ref{ref-4.1-8} is that it provides
 reasonable non-commutative substitutes for
``regularity'', ``birationality'' and ``crepancy''.  Obviously
regularity correspond to the condition $\gldim {{A}} <\infty$.

Let $K$ be the function field of $R$. To have a substitute for birationality 
we note that in non-commutative geometry it is customary to 
replace a
ring with its module category. So birationality should be expressed by the
fact that ${{A}}\otimes_R K$ is Morita equivalent to $K$. I.e. we
should have ${{A}}\otimes_R K=M_n(K)$.

Assume that $R$ is complete local of dimension $n$ and let $\omega_R$
be the dualizing module of $R$. By the Gorenstein hypotheses
$\omega_R$ is an invertible $R$-module. The Grothendieck dualizing
complex  of $R$ is given by $D_R=\omega_R[n]$. By the
adjunction formula ${{A}}$ has a dualizing complex given by
$D_A=\RHom({{A}},D_R)$. It follows that ${{A}}$ has a dualizing
complex concentrated in one degree if an only if $\Ext^i_R({{A}},
\omega_R)=0$ for $i>0$, i.e. if and only if ${{A}}$ is
Cohen-Macaulay.  If ${{A}}$ is Cohen-Macaulay then in particular it
is reflexive.

If ${{A}}$ is Cohen-Macaulay and has finite global dimension then it
is homologically homogeneous. Hence it follows from \cite{BH} that
if $\pi$ is a height one prime in $R$ then ${{A}}_\pi$ is a 
hereditary order in $M_n(K)$ over the discrete valuation ring by $R_\pi$. 
We denote 
the ramification index \cite{reiner} of ${{A}}_\pi$ by $e(\pi)$. 

If ${{A}}$ is Cohen-Macaulay then the dualizing module of ${{A}}$ is
given by $\omega_{{A}}=\Hom_R({{A}},\omega_R)$ and it is easy to prove
that the latter is equal to the tensor product of $\omega_R$ with
$(\otimes_\pi \pi^{-1+e(\pi)})^{\ast\ast}$ where the product runs over
all height one primes in $R$. It follows that that $\omega_{{A}}$ will
be generated by $\omega_R$ (a substitute for crepancy) if and only if
$e(\pi)=1$ for all $\pi$, i.e. if and only if ${{A}}$ is a maximal order.
According to \cite{AG1} all maximal orders in $M_n(K)$ are of the form
$\End_R(M)$ for a reflexive $R$-module $M$. This finishes our
motivating discussion.
\begin{remark} 
\label{ref-4.3-9}
It is tempting to weaken Definition \ref{ref-4.1-8} in such a way
as to require only that ${{A}}$ is an unramified maximal order over $R$.
The proof Theorem \ref{ref-6.3.1-25} below does not work in this added generality.
\end{remark}
\begin{remark} In all our examples the $R$-module $M$ may be taken to be
itself Cohen-Macaulay but we have no proof that this is always possible. Note
that the implication $\End_R(M)$ Cohen-Macaulay $\Rightarrow$ $M$ 
Cohen-Macaulay is false. 
\end{remark}
\begin{remark}
If our basescheme $X$ is not affine then we may define a non-commutative
crepant resolution of $X$ as a stack of abelian categories which is, locally
on an affine open  $\Spec R$, of the form $\Mod({{A}})$ where ${{A}}$ is
a non-commutative crepant resolution of $R$.
\end{remark}
The following conjecture is inspired by the conjecture of Bondal and Orlov.
\begin{conjecture}
\label{ref-4.6-10}
 All crepant resolutions of $X$ (commutative as well as non-commutative
ones) are derived equivalent.
\end{conjecture} 
In subsequent sections we will give some examples of non-commutative
crepant resolutions and in addition we will give some evidence for
Conjecture \ref{ref-4.6-10}. In particular we will prove
Conjecture \ref{ref-4.6-10} for three dimensional
terminal Gorenstein singularities. 
\section{Resolutions with fibers of dimension $\le 1$}
\label{ref-5-11}
The following result was proved is \cite{VdB31} (our hypotheses 
here are slightly more general but this does not affect the proof).
\begin{theorem}
\label{ref-5.1-12} Assume 
that there exists a crepant resolution of singularities  $f:Y\r X$  such that  
the dimensions 
of the fibers of $f$ are $\le 1$ and such that the exceptional locus of $f$ has
codimension $\ge 2$. Then $R$ has a non-commutative crepant
resolution ${{A}}=\End_R(M)$ where $M$ is in addition Cohen-Macaulay.
Furthermore $Y$ and ${{A}}$ are derived equivalent.
\end{theorem}
Let us briefly recall how $M$ is constructed. The Grauert-Riemenschneider
theorem implies $H^i(Y,\Oscr_Y)=0$ for $i>0$ (thus $X$ has rational
singularities).

Let $\Lscr$ be an ample
line bundle on $Y$ generated by global sections. The hypotheses on
the fibers of $f$ imply that $(\Oscr_Y\oplus \Lscr)^\perp=0$ in $D(\Qch(Y))$.
Take an extension
\[
0\r \Oscr_Y^{r}\r \Mscr'\r \Lscr\r 0
\]
associated to a set of $r$ generators of $\Ext^1_Y(\Lscr,\Oscr_Y)$ as $R$-module
 and put
$\Mscr=\Mscr'\oplus\Oscr_Y$. Then $\Mscr^\perp=(\Oscr_Y\oplus \Lscr)^\perp=0$
and furthermore $\Ext^i_A(\Mscr,\Mscr)=0$ for $i>0$.

Put $M=\Gamma(Y,\Mscr)$. The hypotheses imply that ${{A}}=\End_R(M)=
\End_Y(\Mscr)$. Hence $Y$ is derived equivalent to ${{A}}$. In particular
${{A}}$ has finite global dimension. For the fact that ${{A}}$ and
$M$ are Cohen-Macaulay we refer to \cite{VdB31}.

\section{Construction a crepant resolution starting from a 
non-commutative one}
\label{ref-6-13}
\subsection{Introduction}
\label{ref-6.1-14}
In this section we assume that ${{A}}=\End_R(M)$ is an arbitrary
non-commutative crepant resolution of $R$.  We will show that the
beautiful  approach in \cite{BKR} to the McKay
correspondence generalizes to this situation. In fact an almost
literal copy of the proof works, but we will nevertheless give a
summary of it, in order to convince the reader (and
ourselves!) that nothing specific to the situation of $G$-equivariant
sheaves is used in \cite{BKR} (see \S\ref{ref-6.3-22}-\ref{ref-6.6-35}).
The only part that is specific to our more general situation 
is in \S\ref{ref-6.3-22}.

An obvious first step in our more general situation is to take for $Y$
some type of ${{A}}$-Hilbert-scheme of $M$. This would work provided
that $M$ is Cohen-Macaulay.  Unfortunately we don't know if this is
always the case (but as already said, it is true in all examples we
know).

Therefore we take a slightly different approach. We construct $Y$ as 
a moduli-space of certain stable ${{A}}$-representations. It
is standard how to do this but since our base ring $R$ is somewhat more
general than usual we
recall the necessary steps in the 
next section.

\subsection{Moduli spaces of representations}
\label{ref-6.2-15}
In this section $R$ is a commutative noetherian $k$-algebra where as
usual $k$ is algebraically closed of characteristic zero.  Put
$X=\Spec R$.  Let ${{A}}$ be an $R$-algebra which is finitely
generated as $R$-module. Let $(e_i)_{i=1,\ldots,p}$ be pairwise
orthogonal idempotents in ${{A}}$ such that $1=\sum_i e_i$ and let
$D=\oplus_i Re_i\subset {{A}}$ be the corresponding diagonal
subalgebra.

For a map $R\r K$ with $K$ a field and $V$ a finite dimensional
${{A}}\otimes_R K$ representation we write $\udim V=(\dim_K Ve_i)_i\in \ZZ^p$.
We put on $\ZZ^p$ the ordinary Cartesian scalar product.

Pick $\lambda\in \ZZ^p$ and let $\alpha=\udim V$. Following \cite{King}
We say that $V$ is
\emph{(semi-)stable} (with respect to $\lambda$) if  
\begin{equation}
\label{ref-6.1-16}
(\lambda,\alpha)=0
\end{equation}
and if
for any proper
subrepresentation $W$ of $V$ with $\beta=\udim W$ we have
\begin{equation}
\label{ref-6.2-17}
(\lambda,\beta)\,\,(\ge)> 0
\end{equation}
Note that if 
\begin{equation}
\label{ref-6.3-18}
(\lambda,\beta)\neq 0\qquad \text{for $0<\beta<\alpha$}
\end{equation}
then stability and
semi-stability are equivalent. For a fixed $\alpha$ there will exist $\lambda$
satisfying  \eqref{ref-6.1-16} and \eqref{ref-6.2-17} if and
only if the greatest common divisor of all $\alpha_i$ is one. 

An  \emph{affine family} of ${{A}}$ representations
with dimension vector $\alpha$ is a commutative $R$-algebra
$T$ together with a finitely generated ${{A}}\otimes_R T$ module
$P$ which is projective as $T$-module such that $Pe_i$ has
constant rank $\alpha_i$ for all $i$. For such a $P$
we write $\udim P=\alpha$. This is equivalent to saying that
for any map of $T$ to a field $K$ we have  
$\udim(P\otimes_T K)=\alpha$. We say that $P$ is (semi-)stable
if for any $K$ we have that $P\otimes_R K$ is (semi-)stable. 
Non-affine families  are defined in the obvious way by gluing affine
families. We call families \emph{equivalent} if they are locally isomorphic.

Let $\alpha=(\alpha_i)_{i=1,\ldots,p}$  be natural numbers which are relatively
prime and pick a $\lambda$ satisfying  \eqref{ref-6.1-16} and \eqref{ref-6.3-18}.
Consider the functor ${\mathbf{R}}^s$ which assigns to a commutative
$R$-scheme $Z$ the following set:
\begin{multline*}
\{\text{equivalence classes of families of $\lambda$-stable ${{A}}$-representations}\\
\text{over $Z$
with dimension vector $\alpha$}\}
\end{multline*}
\begin{propositions}
The functor ${\mathbf{R}}^s$ is representable by a projective
scheme over~$X$.
\end{propositions}
We need the following lemma.
\begin{lemmas}
  Assume that $R'$ is a commutative $R$-algebra and that $T$ is a
  commutative $R'$-algebra. Put ${{A}}'={{A}}\otimes_R R'$.
 Assume that $P$ is a
  $T$-family of stable ${{A}}'$-representations.  Then 
$P$ is also stable as a a family of
  ${{A}}$-representations.
\end{lemmas}
\begin{proof} This is trivial from the definition.
\end{proof}
If we temporarily write ${\mathbf{R}}^s_R$ for ${\mathbf{R}}^s$ then
this lemma implies that ${\mathbf{R}}^s_{R'}={\mathbf{R}}^s_{R}\times_
{\Spec R}\Spec R'$. Hence if ${\mathbf{R}}^s_{R}$ is representable by a projective
scheme over $R$
then so is ${\mathbf{R}}^s_{R'}$. 

We use this as follows: we may find a finitely generated subring $R_0$
of $R$ and a module finite $R_0$ algebra ${{A}}_0$ such that
${{A}}={{A}}_0\otimes_{R_0} R$. It is then sufficient to prove that
${\mathbf{R}}^s_{R_0}$ is representable by a projective scheme over $R_0$.
Hence without loss of generality we may (and we will) assume that $R$
is finitely generated over $k$.

Put $\bar{\alpha}=\sum_i\alpha_i$.
Define ${{A}}'$ as the centralizer of $M_{\bar{\alpha}}(k)$ in 
$M_{\bar{\alpha}}(R)\ast_{D} {{A}}$.
So we have
\begin{equation}
\label{ref-6.4-19}
M_{\bar{\alpha}}(R)\ast_D {{A}}=M_{\bar{\alpha}}(k)\otimes {{A}}'
\end{equation}
where the obvious copy of 
$M_{\bar{\alpha}}(k)$ on the left and right is the same.

Note that ${{A}}'$ is an $R$-algebra. Since
$M_{\bar{\alpha}}(R)\ast_{D} {{A}}= M_{\bar{\alpha}}({{A}}')$ is
finitely generated the same is true for ${{A}}'$.

Put $S'={{A}}'/[{{A}}',{{A}}']$ and $W'=\Spec S'$. $W'$ is an
affine $R$-scheme of finite type representing the functor \cite{schofield1} 
which
assigns to a commutative $R$-algebra $T$:
\begin{multline}
\label{ref-6.5-20}
\{\text{${{A}}$-module structures on $T^{\bar{\alpha}}$ which commute with the
} \\ \text{$T$-structure such that $e_i$ acts as
$\diag(0^{\alpha_1+\cdots+\alpha_{i-1}},1^{\alpha_i},
0^{\alpha_{i+1}+\cdots+\alpha_p})$}\}
\end{multline}
The corresponding universal bundle is given by $U'_0=k^{\bar{\alpha}}
\otimes_k S$ with the 
right action of ${{A}}$ obtained via the composition ${{A}}\r
M_{\bar{\alpha}}(R)\ast_D {{A}}=M_{\bar{\alpha}}({{A}}')\r M_
{\bar{\alpha}}(S')$.

Let $G=\prod_i \GL_{\alpha_i}(k)$
and let $PG$ be equal to $G$ modulo its center. 
Conjugation on the first factor 
induces a $G$ action on  $M_{\bar{\alpha}}(R)
\ast_D {{A}}$, leaving the elements of
${{A}}$ invariant. 

The  $G$ action on $M_{\bar{\alpha}}(R)\ast_D {{A}}$ leaves 
$M_{\bar{\alpha}}(k)$ stable
and hence it induces rational $G$-actions on ${{A}}'$, $S'$ and 
$W'$. The $G$-action on the righthand side of \eqref{ref-6.4-19}
is the diagonal one.

Putting $g(v\otimes\gamma)=vg^{-1}\otimes g\gamma$ defines a rational
$G$-action on $U'_0$. The center of $G$ acts trivially on $S'$ and via
scalar multiplication on $U'_0$. Let $a\in\ZZ^p$ be such that
$(\alpha,a)=1$ and define $U'=U'_0\otimes_S \bigotimes_i \wedge^{\alpha_i} (U'_0
e_i)^{\otimes -a_i}$.
$U'$ is still a ${{A}}\otimes_R S'$ module which is projective as
$S'$-module and which  has dimension vector $\alpha$.  The center of $G$
now acts trivially on $U'$. 

As usual the stability condition \eqref{ref-6.2-17} corresponds
to a character $\chi:G\r k^\ast$ \cite{King}.  This character defines
a $\ZZ$-grading on $S'$.

Let $G_0=\ker \chi$ and put $S=(S')^{G_0}$, $U=(U')^G$. Then $S$ is still
$\ZZ$-graded and $U$ is a graded $S$-module. Put $W=\Proj S_{\ge 0}$ and let
$f:W\r X$ be the structure map. If $W^{\prime s}$ is the open subset  of $W$
corresponding to the complement of the closed subscheme defined by $S_{>0}$
then it follows from the Luna slice theorem \cite{Luna} that $W^{\prime s}\r W$ is
a principal $PG$-fiber  bundle. 

Let $\Uscr$ be the coherent sheaf on $W$ which corresponds to $U$.
$\Uscr$ is a sheaf of right ${{A}}$-modules. It follows from standard
descent theory that the pullback of $\Uscr$ to $W^{\prime s}$ is the
restriction of $U'$ (considered as a coherent sheaf on $W'$). In particular
$\Uscr$ is a vector bundle on $W$.

It is now standard that $W$ represents the functor ${\mathbf{R}}^s$ and that $\Uscr$ is
the corresponding universal bundle. See \cite{Newstead} for the case of vector 
bundles.
\begin{lemmas} \label{ref-6.2.3-21} If $A=M_{\bar{\alpha}}(R)$   then the
  map $W\r X$ is an isomorphism.
\end{lemmas}
\begin{proof} It follows from Morita theory  that in this case the functor 
 $\mathbf{R}^s$ is represented by $X$. So this proves what we want. 
\end{proof}
\subsection{Application to our situation}
\label{ref-6.3-22}
Now we return to our standard assumption that $R$ is a normal Gorenstein
domain and we let ${{A}}=\End_R(M)$ be an arbitrary
non-commutative crepant resolution of $R$.

Since ${{A}}$ has finite global dimension there is a ${{A}}$-resolution
\begin{equation}
\label{ref-6.6-23}
0\r P\r {{A}}^a\r \cdots \r {{A}}^b\r M\r 0
\end{equation}
with $P$ projective. By the reflexive Morita correspondence we have
$P=\Hom_R(M,M_1)$ where 
\begin{equation}
\label{ref-6.7-24}
M_1\oplus M_2\cong M^{\oplus c}
\end{equation}
for some $c\in\NN$ 
and some other reflexive $R$-module $M_2$.

From \eqref{ref-6.6-23} we obtain $\rk M=d\rk {{A}}\pm \rk P$ for some 
$d\in\ZZ$, where
``$\rk U$'' denotes the rank of $U$ as  $R$-module. Simplifying we find
$1=d\rk M\pm \rk M_1$ and hence the ranks
of $M$ and $M_1$ are coprime. Increasing $c$ if necessary we may also
assume that  the ranks of $M_1$ and $M_2$ are coprime. 

We now replace
$M$ by $M^{\oplus c}$. This changes ${{A}}$ into something Morita
equivalent and the decomposition \eqref{ref-6.7-24} becomes 
$M=M_1\oplus M_2$.

Let us more generally consider a decomposition $M=\oplus_{i=1}^p M_i$.
Denote the rank of $M_i$ by $\alpha_i$ and put $\alpha=(\alpha_i)_i$.
Take $\lambda$ satisfying \eqref{ref-6.1-16} and \eqref{ref-6.3-18} and
let $f:W\r X$ and $\Uscr$ be as in the previous section.

If we let $X_1\subset X=\Spec R$ be the locus where $M$ is locally
free then it follows from lemma \ref{ref-6.2.3-21} that $f^{-1}(X_1)\r X_1$ is
an isomorphism.
We define $Y$ as the unique irreducible component of $W$ mapping onto $X$.
We will denote the restriction of $f$ to $Y$ also by $f$. 
We let $\Mscr$ be the restriction of $\Uscr$ to $Y$. $\Mscr$ is
a sheaf of ${{A}}$-modules on $Y$.

Following \cite{BKR} we now define a pair of adjoint functors between
$D^b(\coh(Y))$ and $D^b(\mod({{A}}))$. 
\begin{gather*}
\Phi:D^b(\coh(Y))\r D^b(\mod({{A}})):C\mapsto 
R\Gamma(C\Lotimes_{\Oscr_Y}\Mscr)\\
\Psi:D^b(\mod({{A}}))\r D^b(\coh(Y)):D\mapsto D\Lotimes_{{A}} \Mscr^\ast
\end{gather*}

The following is a straightforward generalization of \cite{BKR}.
\begin{theorems}
\label{ref-6.3.1-25}
Assume that  for any point $x\in X$ of codimension $n$
the fiber product $(Y\times_X Y)\times_X \Spec \Oscr_{X,x}$ has dimension $\le n+1$. Then
$f:Y\r X$ is a crepant resolution of $X$ and $\Phi$ and $\Psi$ are inverse
equivalences. 
\end{theorems}

To give the proof we need to say something about spanning classes
and Serre functors.  This is done in the next section.
\subsection{Relative Serre duality}
\label{ref-6.4-26}
Assume that $\Dscr\subset \Cscr$ is a full faithful inclusion
of $k$-linear triangulated categories such that for 
$C\in \Cscr$, $D\in\Dscr$ we
have $\sum_i \dim \Ext^i(C,D)<\infty$, $\sum_i \dim \Ext^i(D,C)<\infty$. 
We say  that an auto-equivalence of triangulated categories $S:\Dscr\r \Dscr$
is a relative Serre functor for the pair $(\Dscr,\Cscr)$ if the
following properties hold:
\begin{enumerate}
\item $S$ leaves $\Dscr$ stable.
\item For $C\in \Cscr$, $D\in\Dscr$ there are natural isomorphisms
\begin{equation}
\eta_{D,C}:\Hom(D,C)\r \Hom(C,SD)^\ast
\end{equation}
\end{enumerate}
Thus in particular $S$ is Serre functor for $\Dscr$. 

We will use this in the following situation.
\begin{lemmas} 
\label{ref-6.4.1-27}
Assume that $X=\Spec R$ where $R$ is a complete local 
  noetherian ring and let $f:Y\r X$ be a projective map. Let $\Ascr$
  be a coherent sheaf of $\Oscr_Y$ algebras such that for every $y\in
  Y$ the stalk of $\Ascr$ at $y$ has finite global dimension.  
 Let $\Cscr$ be the bounded derived
  category of $\coh(\Ascr)$ and let $\Dscr$ be the full subcategory of
  complexes whose homology has support in $Y_0=f^{-1}(x)$ where $x\in
  X$ is the closed point. Then the pair $(\Dscr,\Cscr)$ has a relative
  Serre functor given by tensoring with the Grothendieck dualizing
  complex  $D_\Ascr=\uRHom_{Y}(\Ascr,D_Y)$ of $\Ascr$.
\end{lemmas}
\begin{proof}
We claim first that if $D\in \Dscr$, $C\in \Cscr$ then we have the
following identity:
\begin{equation}
\label{ref-6.9-28}
\Hom(D,D_Y)\cong \Gamma(Y,D)^\ast
\end{equation}
Let $Y_n$ be the $n$'th formal neighborhood of $Y_0$ and let
$j_n:Y_n\r X$ be the inclusion map. 

 It is easy to see that there is some $n$
and some $D_n\in D^b(\coh(Y_n))$ such that $D=Rj_{n\ast} D_n$. Then we
have $\Hom(D,D_Y)=\Hom(Rj_{n\ast} D_n,D_Y)=\Hom(D_n,j_n^!(D_Y))=
\Hom(D_n,D_{Y_n})$. Now since $Y_n$ is proper it satisfies classical
Serre duality. I.e. there are isomorphisms $\Hom(D_n,D_{Y_n})\cong
\Gamma(Y_n,D_n)^\ast=\Gamma(Y,D)^\ast$. It is easy to see that
the resulting isomorphism $\Hom(D,D_Y)\cong \Gamma(Y,D)^\ast$ is independent
of $n$. This finishes the proof of \eqref{ref-6.9-28}.

We compute for $D\in \Dscr$, $C\in \Cscr$
\begin{equation}
\label{ref-6.10-29}
\begin{split}
\Hom_\Ascr(D,C)&=\Gamma(Y,\uRHom_\Ascr(D,C))
\\
&\cong\Hom(\uRHom_\Ascr(D,C),D_Y)^\ast\qquad\text{(by \eqref{ref-6.9-28})}
\end{split}
\end{equation}

We claim that there are natural isomorphism 
\begin{equation}
\label{ref-6.11-30}
\uRHom(\uRHom_\Ascr(D,C),D_Y))\cong \uRHom_\Ascr(C,D\Lotimes_\Ascr D_\Ascr)
\end{equation}
This is easily proved by replacing $D_Y$ with a bounded injective
complex and $D$ with a  locally projective complex (using the projectivity
of $f$).

So combining \eqref{ref-6.10-29} and \eqref{ref-6.11-30} we obtain an isomorphism
$\Hom_\Ascr(D,C)\cong \Hom_\Ascr(C,D\Lotimes_\Ascr D_\Ascr)^\ast$ which
finishes the proof.
\end{proof}
\subsection{Spanning classes}
\label{ref-6.5-31}
If $\Cscr$ is a triangulated category then a \emph{spanning
class} \cite{Bridgeland}
 $\Omega\subset \Cscr$ is a set of objects such that $\Omega^\perp=0$
and ${}^\perp\Omega=0$. 
\begin{examples}
\label{ref-6.5.1-32}
Let $f:Y\r X$ be a \emph{proper} map where $X$ is the spectrum of a 
noetherian local
ring with algebraically closed residue field $k$. Then $\Omega=\{\Oscr_y\mid
y\in Y(k)\}$ is a spanning class for $D^b(\coh(Y))$.
\end{examples}
Assume that $F:\Cscr\r \Escr$  is a functor between
pairs of triangulated categories. Assume that $F$ has a left adjoint
$G$ and a right adjoint $H$ and that $\Omega\subset \Cscr$ is a spanning class.
\begin{lemmas} \cite{BO}\cite{Bridgeland}
$F:\Cscr\r\Escr$ is fully faithful if and only if the map
\begin{equation}
\label{ref-6.12-33}
\Ext^i_\Cscr(\omega,\omega')\r \Ext^i_\Escr(F\omega,F\omega')
\end{equation}
is an isomorphism for all $\omega,\omega'\in \Omega$.
\end{lemmas}
Suppose now that $F$ is actually a compatible pair of fully faithful functors
 $(\Dscr,\Cscr)\r (\Fscr,\Escr)$ between pairs of triangulated categories
satisfying the hypotheses given in the beginning of \S\ref{ref-6.4-26}
 which in addition have relative Serre functors $S_\Cscr$ and 
$S_\Escr$. Assume $\Omega\subset \Dscr$. 
\begin{lemmas}
\label{ref-6.5.3-34}
Assume that $\Cscr$ is not trivial and that $\Escr$ is connected. 
Assume in addition that $S_\Escr F\omega=FS_\Cscr\omega$  
for $\omega\in \Omega$.
Then $F$ is an equivalence of categories.
\end{lemmas}
\begin{proof} The proof is the same as in \cite{Bridgeland} with Serre functors
being replaced by relative ones. 
\end{proof}
\subsection{Proof of the main theorem}
\label{ref-6.6-35}
We now complete the proof of Theorem \ref{ref-6.3.1-25}.
We need to show that the canonical natural transformations
$\Psi\Phi\r \Id$ and $\Id\r \Phi\Psi$ are isomorphisms. Since
everything is compatible with base change we may
 reduce to the case where $R$ is a complete local ring
containing a copy of its algebraically closed residue field. We will
denote this new residue field also by $k$. We let $x$ be the unique
closed point of $X=\Spec R$.

The functors $\Phi$ and $\Psi$ have versions for left modules, denoted
by $\Phi^\circ$ and $\Psi^\circ$ respectively,  which are given by
the formulas $R\Gamma(Y,\Mscr^\ast\Lotimes_{\Oscr_Y} -)$ and
$\Mscr\Lotimes_{{A}}-$. It it easy to see that
$\Phi^\circ=\DD_{{A}}\circ \Phi\circ \DD_Y$.  Hence $\Phi$ also has a
right adjoint given by $\DD_Y\circ\Psi^\circ\circ \DD_{{A}}$.
 Furthermore as in Example
\ref{ref-6.5.1-32} the objects $\Oscr_y$ form a spanning class for
$D^b(\coh(Y))$.

For $y\in Y(k)$ denote by $\Mscr_y=\Phi(\Oscr_y)$ 
the fiber of $\Mscr$ at $y$. To prove that $\Phi$ is fully faithful we
need to prove that the canonical maps
\begin{equation}
\label{ref-6.13-36}
\Ext^i_Y(\Oscr_y,\Oscr_{y'})\r \Ext^i_{{A}}(\Mscr_y,\Mscr_{y'})
\end{equation}
are isomorphisms for $y,y'\in Y$.

What do we know already?
\begin{enumerate}
\item \eqref{ref-6.13-36} is certainly an isomorphism for $i=0$.
\item By Serre duality for ${{A}}$ and the fact that
  $D_{{A}}\cong{{A}}[n]$ (lemma \eqref{ref-6.4.1-27})
  \eqref{ref-6.13-36} is an isomorphism if $i=n$ and $y\neq y'$.
\end{enumerate}
There is one more subtle piece of information that may be obtained. 
Recall that $Y$ is a closed subscheme
of the scheme $W$ representing the functor of  stable
${{A}}$-representations.  This
means that there is an injection (the Kodaira-Spencer map):
\[
\phi:\Ext^1_Y(\Oscr_y,\Oscr_y)=T_{Y,y}\hookrightarrow
T_{W,y}=\Ext^1_W(\Mscr_y,\Mscr_y)
\]
where $T_{\ast,y}$ denotes the tangent space at $y$. The map $\phi$ is
constructed as follows. An element of $T_{Y,y}$ corresponds to a map
$u:\Spec k[\epsilon]/(\epsilon^2)\r Y$ and hence to an extension $E$ of
$\Oscr_y$ with itself. Then $\phi(u)=u^\ast(E)$. Thus the Kodaira-Spencer
map coincides with \eqref{ref-6.13-36} for $i=1$ and $y=y'$. Hence
we have 
\begin{enumerate}
\setcounter{enumi}{2}
\item \eqref{ref-6.13-36} is an injection for $i=1$ and $y=y'$.
\end{enumerate}
Using an amazing trick based on the intersection theorem in commutative
algebra it is shown in \cite{BKR} that (1),(2) and (3) are sufficient to
prove that $\Phi$ is fully faithful 
(under the standing hypothesis $\dim Y\times_X Y\le
n+1$).

We need a succinct description of $\Psi\Phi$. Let $Y\hat{\times} Y=
(Y\times Y)_{X\times X}
\Spec \hat{\Oscr}_{X\times X,(x,x)}$. This a noetherian scheme proper
over $\Spec \hat{\Oscr}_{X\times X,(x,x)}$. We denote the projections
$Y\hat{\times} Y\r Y$ by $\pr_{1,2}$. 

$Y\times_X Y$ may be considered as a closed subscheme of $Y\hat{\times}
Y$ (no need to complete).  We define $\Mscr\boxtimes_{{A}}\Mscr^\ast$
as the coherent sheaf on $Y\hat{\times} Y$ such that for affine opens
$U,V\subset Y$ we have
$(\Mscr\boxtimes_{{A}}\Mscr^\ast)(U\hat{\times} V)=
\Mscr(U)\otimes_{{A}} \Mscr(V)^\ast$. $\Mscr\boxtimes_{{A}}\Mscr^\ast$ is
clearly supported on $Y\times_X Y$. Using suitable flat resolutions we
may define the analogous derived object
$\Qscr=\Mscr\Lboxtimes_{{A}}\Mscr^\ast$ which is also supported on
$Y\times_X Y$. Then it is easy to see that 
\[
\Psi\Phi=R\pr_{2\ast}(L\pr_1^{\ast}(-)\otimes_{Y\hat{\times} Y}\Qscr)
\]
We have $\Oscr_{y,y'}\otimes_{Y\hat{\times} Y}
\Qscr=\Mscr_y\Lotimes_{{A}} \Mscr^\ast_{y'}$. Using lemma
\ref{ref-6.4.1-27} for ${{A}}$ we have $\Mscr^\ast_{y'}=
\Hom_{{A}}({{A}},\Mscr_{y'})^\ast= \Hom_{{A}}(\Mscr_{y'},{{A}}[n])$
and since the righthand side of the last equality has homology only in
degree zero, it is equal to $\RHom(\Mscr_{y'},{{A}}[n])$. Thus
$\Mscr_y\Lotimes_{{{A}}}
\Mscr^\ast_{y'}=\RHom_{{A}}(\Mscr_{y'},\Mscr_y)[n]$.

Thus if $y\neq y'$ then $H^{-n}(\Oscr_{y,y'}
\Lotimes_{Y\hat{\times} Y}
 \Qscr)=0$ and using relative Serre duality for ${{A}}$ again we also have
$H^0(\Oscr_{y,y'}
\Lotimes_{Y\hat{\times} Y}
 \Qscr)=\Ext^n_{{A}}(\Mscr_y,\Mscr_{y'})=\Hom_{{A}}(\Mscr_{y'},\Mscr_{y})=0$.
 So the range of possible non-zero values for $H^i(\Oscr_{y,y'}
\Lotimes_{Y\hat{\times} Y}
 \Qscr)=0$ has size $n-1$. But this implies by
 the intersection theorem (see \cite{BKR,BrIy}) that on the complement of the
diagonal $\Qscr$ has support of dimension $\ge \dim(Y\hat{\times} Y)+1-(n-1)
=n+2$ (if non-empty). But by hypotheses the support
of $\Qscr$ has dimension less than or equal the dimension of $Y\times_X Y$
which is $n+1$.
We conclude
that $\Qscr$ is supported on $\Delta$. 

Consider $\Psi\Phi\Oscr_y$. By the previous paragraph this is a
complex supported on $y$ living in non-positive degree and we have
$\Ext^i(\Psi\Phi\Oscr_y,\Oscr_y)= \Ext^i(\Mscr_y,\Mscr_y)$ which is
non-zero only in degrees $[0,\ldots,n]$.  We claim
$H^0(\Psi\Phi\Oscr_y)=\Oscr_y$.

Let $c_y[1]$ be the cone over $\Psi\Phi\Oscr_y\r \Oscr_y$. Thus we have
a triangle
\begin{equation}
\label{ref-6.14-37}
c_y\r \Psi\Phi\Oscr_y\r \Oscr_y
\end{equation}
 By the previous
paragraph $c_y$ is supported in $y$.
If is also easy to see that the homology
of $c_y$ is concentrated in non-positive degree.   There is an exact sequence
\[
0\r \Hom(\Oscr_y,\Oscr_y)\r \Hom(\Mscr_y,\Mscr_y)\r \Hom(c_y,\Oscr_y)\r
\Ext^1(\Oscr_y,\Oscr_y)\r \Ext^1(\Mscr_y,\Mscr_y)
\]
and using (1)(3) we conclude $\Hom(c_y,\Oscr_y)=0$ and hence
$H^0(c_y)=0$.  The fact that $H^0(\Psi\Phi\Oscr_y)=\Oscr_y$ now
follows from \eqref{ref-6.14-37}. By the intersection theorem (see \cite{BKR,BrIy})
we conclude that $Y$ is regular at $y$ and that $\Psi\Phi\Oscr_y=\Oscr_y$.

Since this is true for all $y$ we now know that $Y$ is regular and that
\eqref{ref-6.13-36} holds. Thus $\Phi$ is faithful. 

To prove that $\Phi$ is an equivalence we use lemma \ref{ref-6.5.3-34} since
$D^b(\mod({{A}}))$ is
trivially connected. 
We need that $\Phi S_Y \Oscr_y\cong S_{{A}}  \Mscr_y$. Now $S_{{A}}$ is just
shifting $n$ places to the left, and since $Y$ is regular
 $\omega_Y$ is invertible and thus
 $S_Y\Oscr_y$ is just $\Oscr_y[n]$.
 
 Finally we prove that $f$ is crepant.  It is sufficient to prove
that $\omega_Y\cong \Oscr_Y$. Indeed if this is the case then $
f_\ast\omega_Y\cong \Oscr_X$ is reflexive and hence it is equal to $\omega_X$.
Furthermore it is then also clear that $f^\ast\omega_X=\omega_Y$.

Let $D^b_x(\coh(Y))$ the full
subcategory of $D^b(\coh(Y))$ consisting of complexes supported in $f^{-1}(x)$.
Similarly let $D^b_x(\mod({{A}}))\subset D^b(\mod({{A}}))$ be the
complexes supported on $x$. The functor $\Phi$ and $\Psi$ define inverse
equivalences between $D^b_x(\coh(Y))$ and $D^b_x(\mod({{A}}))$.

On $D^b_x(\coh(Y))$ we have that $S_Y[-n]$ is isomorphic to the
identity functor since the same holds for $D^b_x(\mod({{A}}))$. Thus
if $Y_0=f^{-1}(x)$ then $\omega_Y/\omega_Y(-nY_0)\cong 
\Oscr_Y/\Oscr_Y(-nY_0)$. Hence if $\hat{Y}$ is the formal scheme associated
to $Y_0$ then $\hat{\omega}_Y\cong \hat{\Oscr}_Y$. But by the
Grothendieck existence theorem this implies $\omega_Y\cong\Oscr_Y$. 
\begin{remarks}
  As noted in \cite{BKR} the fact that $Y$ is smooth and the fact that
  for $y\in Y$ we have $T_{Y,y}=\Ext^1_Y(\Oscr_y,\Oscr_y)
  \cong\Ext^1_{{A}}(\Mscr_y,\Mscr_y)=T_{W,y}$ implies that $Y$ is a
  connected component of $W$. If $\dim X=3$ then it is shown that
  actually $W=Y$. This result generalizes probably not to our current
  situation.  However if $R$ is complete local and we take a
  decomposition of $M$ into indecomposables (unique up
to non-unique isomorphism by the
  Krull-Schmidt-theorem) then the proof goes through virtually
  unmodified.
\end{remarks}
\begin{remarks} Varying $\lambda$ we get many different crepant 
resolutions of $X$. They are all derived equivalent since they
are all derived equivalent to ${{A}}$. This gives another instance
where the Bondal-Orlov conjecture is true.
\end{remarks}
Now we may prove.
\begin{theorems}
\label{ref-6.6.3-38}
 Assume that $R$ is three-dimensional and has terminal
singularities. 
\begin{enumerate}
\item $R$ has a non-commutative crepant resolution 
if and only if it has a commutative one.
\item Conjecture \ref{ref-4.6-10} is true in this case.
\end{enumerate}
\end{theorems}
\begin{proof}
\begin{enumerate}
\item This follows from Theorems \ref{ref-5.1-12} and \ref{ref-6.3.1-25}.
\item
Let $Y\r X$ be a crepant resolution of singularities and let ${{A}}$
be a non-commutative one.   By Theorem \ref{ref-6.3.1-25} there is another
crepant resolution $Y'\r X$ of $X$ associated to ${{A}}$. The resulting
birational map $Y-\r Y'$ is a composition of flops. Hence by \cite{Br1} 
$Y$ and $Y'$ are derived equivalent. Since $Y'$ and ${{A}}$ are also
derived equivalent, we are done. \qed
\end{enumerate}
\def\qed{}\end{proof}

\section{Cones over del Pezzo surfaces}
\label{ref-7-39}
A standard example of a canonical singularity which is not terminal is
the cubic cone $w^3+x^3+y^3+z^3=0$. This singularity has a crepant
resolution obtained by blowing up the origin. Our aim in this section
is to show that it also has a non-commutative crepant resolution. The
method used applies more generally to cones over Del-Pezzo surfaces
(see Proposition \ref{ref-7.3-45}) below. 

Below $Z$ is regular connected projective scheme of dimension $n-1>0$ with 
an ample line bundle $\Lscr$.
We denote the homogeneous coordinate ring  $\oplus_i  \Gamma(Z,\Lscr^i)$ of 
$Z$  
corresponding to $\Lscr$ 
by $R$. $R$ is a finitely generated
normal graded ring. Put $X=\Spec R$. We will call $R$ a cone over $Z$.
$X$ has unique singularity at
the origin $o$. This singularity has a standard resolution by 
$Y=\underline{\Spec} S(\Lscr)$. Denote the structure map $Y\r Z$ by
$\pi$ and let $f$ be the canonical map $Y\r X$. Thus we have the following
diagram
\[
\begin{CD}
Y @>\pi>> Z\\
@V f VV \\
X
\end{CD}
\]
Let  $X'\subset X$,
$Y'\subset Y$ be the open subsets respectively defined by the ideals 
$R_{>0}$ and 
$(S\Lscr)_{>0}$.  $X'$ is regular and $X-X'=\{o\}$. 
We also have $Y-Y'=E$ where $E$ is the image of
the zero section on $\pi$. The map $f$ restricts to an isomorphism $Y'\r X'$
and so the exceptional locus of $f$ is given by $E$.

We now state a list of properties of $X$ which are well-known and which
are easy to prove. 
\begin{lemma}
\label{wellknown}
\begin{enumerate}
\item $X$ is Cohen-Macaulay if and only if 
\[
H^i(Z,\Lscr^j)=0\qquad\text{\textrm{for all  $0<i<n-1,j\ge 0$}}
\]
\item
$X$ has rational singularities if and only if 
\[
H^i(Z,\Lscr^j)=0\qquad\text{\textrm{for all $i>0,j\ge 0$}}
\]
\item
$X$ has an invertible canonical bundle if and only if $\omega_Z=\Lscr^{-m}$
for $m\in\ZZ$. In that case $\omega_Y=f^\ast\omega_X((m-1)E)$.
\item If $X$ has Gorenstein rational singularities then $\omega_Z=\Lscr^{-m}$
for some $m>0$.
\item If $X$ has Gorenstein rational singularities
then $f$ is not crepant
(or equivalently $X$ is terminal) if and only if $\omega_Z\not\cong\Lscr^{-1}$.
\end{enumerate}
\end{lemma}
Now assume that there is a vector bundle $\Escr_0$ on $Z$ which is
a tilting object.  Put $\Escr=\pi^\ast \Escr_0$. Then $\Escr$ is a tilting
object on $Y$ if and only if 
\begin{equation}
\label{ref-7.1-40}
H^i(Z,\HEnd(\Escr_0,\Escr_0)\otimes\Lscr^j)=0\qquad\text{for $i>0$ and 
$j\ge 0$}
\end{equation}
Thus if \eqref{ref-7.1-40} holds and we put ${{A}}=\End(\Escr)$ then 
${{A}}$ and $Y$ are derived equivalent. 

We claim the following result.
\begin{proposition}  Assume that $\Escr_0$ is a vector bundle on $Z$ which
is a generator of $D(\Qch(Z))$ such that $\Escr=\pi^\ast\Escr_0$ is 
a tilting object on $Y$. Then
$A=\End(\Escr)$ is a non-commutative crepant resolution
of $X$ if and only if the following condition holds:
\begin{equation}
\label{ref-7.2-41}
H^i(Z,\Ascr_0\otimes\Lscr^j)=0\qquad\text{for $i<n-1$ and $j<0$}
\end{equation}
for $\Ascr_0=\HEnd(\Escr_0)$.
\end{proposition}
If $\omega_Z=\Lscr^{-1}$ then \eqref{ref-7.2-41} is always satisfied.
\begin{proof}
  For $\Mscr\in \Qch(Z)$ put $\underline{\Gamma}(Z,\Mscr)=\oplus_{j\in
    \ZZ}\Gamma(Z, \Mscr\otimes\Lscr^j)$ and denote the derived functors of
$\underline{\Gamma}$ by $\underline{H}^\ast$.

If $\Mscr$ is associated to $M\in \Gr(R)$ then it is well-known that we have 
\begin{equation}
\label{ref-7.3-42}
\underline{H}^i(Z,\Mscr)=H^{i+1}_{R_{>0}}(M)\qquad\text{for $i>0$}
\end{equation}
and there is a long exact sequence
\begin{equation}
\label{ref-7.4-43}
0\r H^0_{R_{>0}}(M)\r M\r \underline{\Gamma}(Z,\Mscr)\r H^1_{R_{>0}}(M)\r 0
\end{equation}
Now we prove the first part of the proposition. Since $A$ is derived
equivalent to $Y$
we already know 
$\gldim A<\infty$, 
so we only have to be concerned with the Cohen-Macaulayness of $A$.
The latter is equivalent to $H^i_{R>0}(A)=0$ for $i\le n-1$ (since
the dimension of $R$ is $n$). 

Since we assume \eqref{ref-7.1-40}  the condition \eqref{ref-7.2-41} is
equivalent to 
\begin{equation}
\label{ref-7.5-44}
\begin{split}
\underline{H}^i(Z,\Ascr_0)&=0\qquad\text{for $0<i<n-1$}\\
\underline{\Gamma}(Z,\Ascr_0)_{<0}&=0
\end{split}
\end{equation}
Using \eqref{ref-7.3-42} and \eqref{ref-7.4-43} with $M=A$ we see that $A$ is
Cohen-Macaulay if and only if  $\underline{H}^i(Z,\Ascr_0)=0$ for $0<i<n-1$ and
$A=\underline{\Gamma}(Z,\Ascr_0)$. These conditions correspond
precisely to the conditions given in \eqref{ref-7.5-44}.

Now assume that \eqref{ref-7.5-44} holds and put 
$E=\underline{\Gamma}(\Escr_0)^{\ast\ast}$ where $(-)^{\ast\ast}$ denotes the $R$-bidual. Then $A$ and $\End_R(E)$ are
reflexive $R$-modules which have the same restriction to $X'$ 
(when considered as sheaves). Hence they are equal.

If $\omega_Z=\Lscr^{-1}$ then by  Serre duality we have 
$H^i(Z,\Ascr_0\otimes\Lscr^j)=
H^{n-1-i}(Z,\Ascr_0^\ast \otimes \Lscr^{-j}\otimes\omega_Z)^\ast
=H^{n-1-i}(Z,\Ascr_0\otimes \Lscr^{-j-1})^\ast
$. If $j<0$ and $i<n-1$ then
$-j-1\ge 0$ and $i>0$. I.e. \eqref{ref-7.2-41} follows from \eqref{ref-7.1-40}. 
\end{proof}
Our aim is now to apply this is the case that $Z$ is a surface with 
ample anti-canonical bundle (i.e.\ a Del-Pezzo surface).
Recall that by \cite{Beauville} $Z$ is either $\PP^1\times\PP^1$ or else
is obtained by blowing up $\PP^2$ in $\le 8$ points in general position.
\begin{proposition}
\label{ref-7.3-45} Let $Z$ be a Del-Pezzo surface and let $R$ be a cone
over $Z$ with trivial canonical bundle. Then $R$ has a non-commutative 
crepant resolution.
\end{proposition}
\begin{proof}
Let us first discuss the cases where $\omega_Z$ is a \emph{proper} multiple
of a line bundle. 
If $F$ is an exceptional curve on a surface $Z$ then 
$\deg (\omega_Z\mid F)=-1$ and so $\omega_Z$ cannot be a proper multiple
of a line bundle. The Del-Pezzo surfaces without exceptional curves
are $\PP^2$ and   $\PP^1\times \PP^1$. In the first case 
$\omega_Z=\Oscr_Z(-3)$ and hence  $\Lscr=\Oscr_Z(1)$
and in the second case $\omega_Z=\Oscr_Z(-2,-2)$ and hence
$\Lscr=\Oscr_Z(1,1)$.

In the first case the cone over $Z$ is a 
polynomial ring so this is trivial. In the second case the cone
is given by $R=k[u,v,x,y]/(uv-xy)$ which is standard. The non-commutative
crepant resolution is given by
\[
\begin{pmatrix} R & I \\
I^{-1} & R
\end{pmatrix} 
\]
where $I=(u,x)$.

So from now on we assume $\Lscr=\omega_Z^{-1}$. We will construct a
generator for $D(\Qch(Z))$ satisfying condition \eqref{ref-7.1-40}.  
In fact this follows easily from the results in \cite{KuOr} (I thank
Tom Bridgeland for pointing this out to me). Let $K$ be the canonical
divisor.

It is easy to construct an exceptional collection of vector bundles on
$Z$ generating the derived category $D(\Qch(Z))$. By \cite[Claim
6.5]{KuOr} one may, using a process called ``mutation'',
construct from the initial exceptional collection a new exceptional
collection $E_1,\ldots,E_n$ of vector bundles on $Z$ generating the
derived category $D(\Qch(Z))$ such that the slope function
\[
\mu(E)=-\frac{c_1(E)\cdot K}{r(E)}
\]
is non-decreasing on the sequence
\[
\ldots,E_{n-1}(K),E_n(K),E_n,\ldots,E_n,E_1(-K),E_2(-K),\ldots
\]
and such that any interval of length $n$ is an exceptional collection
of ``type $\Hom$'' (i.e.\ it is an exceptional collection without
higher $\Ext$'s; this is also called a strong exceptional collection).

We claim that the whole sequence has no forward $\Ext^2$'s and no
backward $\Hom$'s and furthermore that all $\Ext^1$'s vanish.  To simplify
the notation we write $E_i(-jK)=E_{i+jn}$.

Let us first consider $\Ext^1$. By Serre duality we have
\[
\dim \Ext^1(E_u,E_v)=\dim \Ext^1(E_v,E_{u-n})
\]
According to \cite[Lemma 3.7]{KuOr} we have
for $v\ge u$ (as $\mu(v)\ge \mu(u)$):
\[
\Ext^1(E_v,E_u)=0\Rightarrow \Ext^1(E_u,E_v)=0
\]
and thus
\[
\Ext^1(E_u,E_{v-n})=0\Rightarrow \Ext^1(E_u,E_v)=0
\]
\[
\Ext^1(E_v,E_u)=0\Rightarrow \Ext^1(E_v,E_{u-n})=0
\]
It is now easy to see that the fact $\Ext^1(E_u,E_v)=0$ for $|u-v|<n$ implies
that all $\Ext^1$'s are zero.

Now we consider $\Hom$'s. We need $\Hom(E_v,E_{u})=0$ for $v>
u$. 
Since $-K$ is effective
on a Del Pezzo surface we have $E_{u-n}\hookrightarrow E_u$. Thus
\begin{equation}
\label{hompropagation}
\Hom(E_v,E_{u})=0\Rightarrow \Hom(E_v,E_{u-n})=0
\end{equation}
We have $\Hom(E_v,E_{u})=0$ for $v-n<u<v$ as $E_{v-n+1},\ldots,E_v$ is
exceptional and also $\Hom(E_v,E_{u})=0$ for $u=v-k$ by Serre duality. Hence
\eqref{hompropagation} implies $\Hom(E_v,E_{u})=0$ for all $u<v$.

Now we consider $\Ext^2$. Again by Serre duality we have 
\[
\dim \Ext^2(E_u,E_v)=\dim \Hom(E_v,E_{u-n})
\]
Hence we need $\Hom(E_v,E_{u-n})=0$ for $v\ge u$. Since $v\ge u$ implies
$v>u-n$ this follows from the vanishing of $\Hom$ as discussed in the previous
paragraph.

\medskip

Put $\Escr_0=E_1\oplus\cdots \oplus E_n$.  The clearly $\Escr_0$
satisfies the conditions of \eqref{ref-7.1-40}.
\end{proof}
\begin{remark} In the published version of this paper \cite{VdB32} we
  considered a more naive tilting object. Assume that $Z$ is obtained
  by blowing up $\PP^2$ in $x_1,\ldots,x_p$ with $p\le 8$ and denote
  the corresponding exceptional curves by $F_1,\ldots,F_p$. Let
  $\alpha:Z\r \PP^2$ be the structure map.  In \cite{VdB32} we took
  $\Escr_0=\alpha^\ast(\Oscr_{\PP^2})\oplus
  \alpha^\ast(\Oscr_{\PP^2}(1))\oplus \alpha^\ast(\Oscr_{\PP^2}(2))
  \oplus \Oscr_Z(F_1) \oplus \cdots \oplus\Oscr_Z(F_p)$ as tilting
  object on $Z$. In contrast to what we stated this tilting object
  does not satisfy the conditions of \eqref{ref-7.1-40}. Indeed 
for $p>4$ we have $H^1(Z,\HEnd(\Escr_0,\Escr_0)\otimes_Z \omega^{-1}_Z)\neq 0$
as one easily checks that $H^1(Z,\operatorname{\mathcal{H}\mathit{om}}(\alpha^\ast(\Oscr_{\PP^2}(2)),
\Oscr_{Z}(F_1)\otimes_Z \omega^{-1}_Z))\neq 0$. 
%
%

\end{remark}

\section{One-dimensional torus invariants}
\label{ref-8-46}
Let $T=k^\ast$ be a one-dimensional torus acting on a finite dimensional vector
space $V$. We may choose a basis $x_1,\ldots,x_n$ for $V$ such that $T$ acts
diagonally: $z\cdot x_i=z^{a_i}x_i$ for some $a_i\in \ZZ$. Put $S=\Sym(V)$ and 
$R=S^T$. In order to avoid trivialities we assume that there are at least 
two strictly positive and two strictly negative $a_i$'s and that the 
greatest common divisor of the $a_i$'s  is one.
Put $N^+=\sum_{a_i>0}a_i$, $N^-=-\sum_{a_i<0}a_i$ and $N=\min(N^+,N^-)$.

It will be convenient to use the Artin-Zhang Proj of a graded ring \cite{AZ}. 
 If
$T$ is a noetherian $\ZZ$-graded ring then $X=\Proj_{\text{AZ}} T$ is
the Grothendieck category $\Gr(T)/\Tors(T)$ where $\Tors(T)$ is the
localizing subcategory of $\Gr(T)$ given by
the graded modules which are limits of right bounded ones. Below
we denote the quotient functor $\Gr(T)\r \Gr(T)/\Tors(T)$ by $\pi$.
We write $\Oscr_X=\pi S$ and $\underline{\Gamma}(X,\Mscr)=
\oplus_{n\in\ZZ}\Hom(\Oscr_X,\Mscr(n))$ where the functor $\Mscr\mapsto
\Mscr(n)$ is obtained from the corresponding functor on $\Gr(T)$. The
derived functors of $\underline{\Gamma}$ are denoted by $\underline{H}^\ast$. 
We write $\coh(X)$ for the category of noetherian objects in $X$.

It is shown
in \cite{AZ} that $\Proj_{AZ} T=\Proj_{AZ} T_{\ge 0}$ so in principle
one may restrict to $\NN$-graded rings but below it will
be more convenient to work with $\ZZ$-graded rings.

Putting $\deg x_i=a_i$ defines a $\ZZ$-grading on $S$ such that
$R=S_0$. By \cite{Hochster} $R$ is Cohen-Macaulay. The other homogeneous
components $S_a$ of $S$ are finitely generated $R$-modules. It is
known \emph{precisely} when they are Cohen-Macaulay \cite{St1,VdB1}.
 In particular one has the following
result.
\begin{lemma} 
\label{ref-8.1-47} $S_a$ is Cohen-Macaulay for $-N^+<a<N^-$.
\end{lemma}
Let $S^+=S$ as graded rings and let $S^-$ be the ring $S$ with modified
grading given by $S^-_n=S_{-n}$ We define $X^\pm=\Proj_{\text{AZ}} S^\pm$. 
Thus if $I^+$, $I^-$ are respectively the ideals generated by $S_{>0}$
and $S_{<0}$ then $X^\pm$ are simply the Grothendieck categories
$\Gr(S^\pm)/\Tors(S^\pm)$ where $\Tors(S^\pm)$ consists of the objects
in $\Gr(S^\pm)$ which are (elementwise) annihilated by powers of
$I^\pm$. Note that $I^+$, $I^-$ are the graded ideals in $S$
respectively generated by $(x_i)_{a_i>0}$ and $(x_i)_{a_i<0}$. From
the description of $\Tors(S^\pm)$ as torsion theories associated to
ideals, it follows that $\Tors(S^\pm)$ is stable.
I.e. $\Tors(S^\pm)$ is  closed under injective hulls \cite{stenstrom}.  In particular
the following lemma follows.
\begin{lemma}
$X^+$, $X^-$ have finite global dimension. 
\end{lemma}
\begin{lemma}
 We have $\underline{H}^i(X^\pm,\Oscr_{X^\pm})=
H^{i+1}_{I^{\pm}}(S^{\pm})$ for $i>0$ and $S^\pm=\underline{\Gamma}(X^{\pm},
\Oscr_{X^\pm})$.
\end{lemma}
\begin{proof} This is a variant on \eqref{ref-7.3-42}\eqref{ref-7.4-43}. Note that 
  $H^{0,1}_{I^\pm}(S^\pm)=0$ because of the hypotheses on the weights
  $(a_i)_i$.
\end{proof}
\begin{lemma} 
\label{ref-8.4-48}
Let $0\le i,j<N^{\pm}$. Then $\Ext^p_{X^\pm}(\Oscr_{X^\pm}(m),
\Oscr_{X^\pm}(n))=0$ for $p>0$ and $\Hom_{X^\pm}(\Oscr_{X^\pm}(m),
\Oscr_{X^\pm}(n))=S^{\pm}_{n-m}$.
\end{lemma}
\begin{proof}
  By the previous lemma we have $\Ext^p_{X^\pm}(\Oscr_{X^\pm}(m),
  \Oscr_{X^\pm}(n))=H^{p+1}_{I^\pm}(S^\pm)_{n-m}$ for $p>0$. A
  standard computation reveals that $H^u_{I^\pm}(S^\pm)$ is zero in
  degrees $> -N^\pm$ which is what we want. 

Also by the previous lemma we find $\Hom_{X^\pm}(\Oscr_{X^\pm}(m),
  \Oscr_{X^\pm}(n))=S^\pm_{n-m}$.
\end{proof}
\begin{lemma}
\label{ref-8.5-49}
 The objects $\Oscr_{X^\pm}(m)$, $m=0,\ldots,N^\pm-1$ form a 
  generating family for $D(X^\pm)$.
\end{lemma}
\begin{proof}
  Let us work with $X^+$.  With a variant of \cite[Lemma
  4.2.2]{BondalVdb} we see that $(\Oscr_{X^+}(n))_n$ generates
  $D(X^+)$. Now we look at the Koszul exact sequence associated to
  $(x_i)_{a_i>0}$.
\[
0\r S(-N^+)\r \cdots \r \oplus_{a_i>0} S(-a_i)\r S\r S/I^+\r 0
\]
Note that $S/I^+$ is right bounded. Shifting and applying $\pi$ we
obtain exact sequences
\[
0\r \Oscr_{X^+}(n-N^+)\r \cdots \r \oplus_{a_i>0}
 \Oscr_{X^+}(n-a_i)\r \Oscr_{X^+}(n)\r 0
\]
Hence it follows that all $\Oscr_{X^+}(n)$ may be obtained using triangles
from $\Oscr_{X^+}(m)$, $m=0,\ldots,N^+-1$ which finishes the proof.
\end{proof}
To simplify the notations below we will now invert, if necessary, the
signs of the $a_i$'s, to insure that $N^+\le N^-$. Thus
$N=N^+$.

Define $\Escr^\pm=\oplus_{m=0,\ldots,N-1}\Oscr_{X^\pm}(m)$, 
 According to lemma \ref{ref-8.4-48}
we have $A\overset{\text{def}}{=}\End_{X^+}(\Escr^+)=\End_{X^-}(\Escr^-)$.

Lemmas \ref{ref-8.4-48},\ref{ref-8.5-49} now yield the following 
\begin{theorem}
  The functor $\RHom(\Escr^+,-)$ defines an equivalence $D(X^+)\r
  D(A)$ and the functor $-\otimes_A \Escr^-$ defines a full
  faithful embedding $D(A)\r D(X^-)$. In particular there is a
full faithful embedding $D(X^+)\r D(X^-)$. All embeddings restrict to
embeddings between the corresponding bounded derived categories of coherent
objects.
 If $N^-=N^+$ (or equivalently
$\sum_i a_i=0$) then all embeddings are equivalences.
\end{theorem}
The existence of an embedding/equivalence $D(X^+)\r D(X^-)$ was proved
 by Kawamata in \cite{Kawamata} (in a slightly more general situation).

We also have
\begin{proposition} $A$ is Cohen-Macaulay.
\end{proposition}
\begin{proof}
According to lemma \ref{ref-8.4-48} we have $A=(S_{n-m})_{0\le m,n< N}$. 
Cohen-Macaulauness now follows from lemma \ref{ref-8.1-47}.
\end{proof}
Now we restrict to the case $N^+=N^-$, or equivalently $\sum_i a_i=0$.
 Under these hypotheses
we have
\begin{lemma} 
\label{ref-8.8-50} $S_0$ is Gorenstein and the $S_a$ are reflexive $S_0$-modules,
satisfying $(S_aS_b)^{\ast\ast}=S_{a+b}$. 
\end{lemma}
\begin{proof} $S_0$ is Gorenstein because of \cite[13.3]{St5}. The hypotheses
that the greatest common divisor of the $(a_i)_i$ is one implies that
  that the generic stabilizer of the $T$-action on
  $W=\Spec \Sym(V)=V^\ast$ is trivial.

In general if $(\zeta_i )_i$ is a point in $W$ then the order
of its stabilizer is equal to the greatest common divisor of the $a_i$'s
such that $\zeta_i\neq 0$. From the relation $\sum_i a_i=0$ we then deduce
that the complement of the locus $W'\subset W$ where $T$ has trivial 
stabilizer, has codimension at least two. The equality $(S_aS_b)^{\ast\ast}=
S_{a+b}$ is now true on $W'$  and hence it is true on $W'\quot T$ since
$W'\r W'\quot T$ is a principal $T$-bundle. Then is is also true on $W/\quot T$
since $W\r W\quot T$ contracts no divisor and hence all modules of covariants
are reflexive \cite[Section 1.3]{Brion}.
\end{proof}
We may now state our final result.
\begin{theorem} Let the notations be as above and assume $\sum_i a_i=0$.
Then $R=S^T$ is a Gorenstein ring with a non-commutative crepant resolution 
given by $A=\End_R(\oplus_{a=0}^{N-1} S_a)$ where $N=\sum_{a_i>0} a_i
=-\sum_{a_i<0} a_i$.
\end{theorem}
\begin{proof} We have $A=\oplus_{m,n=0,\ldots,N-1} (S_{n-m})$ which
  by lemma \ref{ref-8.8-50} is equal to
  $\End_R(\oplus_{a=0}^{N-1} S_a)$.
\end{proof}

\begin{thebibliography}{10}

\bibitem{AZ}
M.~Artin and J.~J. Zhang, {\em Noncommutative projective schemes}, Adv. in
  Math. {\bf 109} (1994), no.~2, 228--287.

\bibitem{Artinx}
M.~Artin, {\em Maximal orders of global dimension and {K}rull dimension two},
  Invent. Math. {\bf 84} (1986), 195--222.

\bibitem{AS}
M.~Artin and W.~Schelter, {\em Graded algebras of global dimension 3}, Adv. in
  Math. {\bf 66} (1987), 171--216.

\bibitem{ATV1}
M.~Artin, J.~Tate, and M.~Van~den Bergh, {\em Some algebras associated to
  automorphisms of elliptic curves}, The Grothendieck Festschrift, vol.~1,
  Birkh\"auser, 1990, pp.~33--85.

\bibitem{ATV2}
\bysame, {\em Modules over regular algebras of dimension 3}, Invent. Math. {\bf
  106} (1991), 335--388.

\bibitem{AG1}
M.~Auslander and O.~Goldman, {\em Maximal orders}, Trans. Amer. Math. Soc. {\bf
  97} (1960), 1--24.

\bibitem{Batyrev}
V.~V. Batyrev, {\em Birational {C}alabi-{Y}au {$n$}-folds have equal {B}etti
  numbers}, New trends in algebraic geometry (Warwick, 1996) (Cambridge),
  London Math. Soc. Lecture Note Ser., vol. 264, Cambridge Univ. Press,
  Cambridge, 1999, pp.~1--11.

\bibitem{Beauville}
A.~Beauville, {\em Surface alg\'ebriques complexes}, Ast{\'e}risque, vol.~54,
  Soc. Math. France, 1978.

\bibitem{BO1}
A.~Bondal and D.~Orlov, {\em Derived categories of coherent sheaves},
  math.AG/0206295.

\bibitem{Bondal1}
\bysame, {\em Semi-orthogonal decompositions for algebraic varieties},
  alg-geom/950601.

\bibitem{BondalVdb}
A.~Bondal and M.~Van~den Bergh, {\em Generators and representability of
  functors in commutative and noncommutative geometry}, Mosc. Math. J. {\bf 3}
  (2003), no.~1, 1--36, 258.

\bibitem{BO}
N.~Bourbaki, {\em Alg\`ebre commutative}, Hermann, Paris, 1960-65.

\bibitem{Bridgeland}
T.~Bridgeland, {\em Equivalences of triangulated categories and
  {F}ourier-{M}ukai transforms}, Bull. London Math. Soc. {\bf 31} (1999),
  no.~1, 25--34.

\bibitem{Br1}
\bysame, {\em Flops and derived categories}, Invent. Math. {\bf 147} (2002),
  613--632.

\bibitem{BKR}
T.~Bridgeland, A.~King, and M.~Reid, {\em The {M}c{K}ay correspondence as an
  equivalence of derived categories}, J. Amer. Math. Soc. {\bf 14} (2001),
  no.~3, 535--554.

\bibitem{BrIy}
T.~Bridgeland and S.~Iyengar, {\em A criterion for regularity of local rings},
  C. R. Math. Acad. Sci. Paris {\bf 342} (2006), no.~10, 723--726.

\bibitem{Brion}
M.~Brion, {\em Sur les modules de covariants}, Ann. Sci. {\'E}cole Norm. Sup.
  (4) {\bf 26} (1993), 1--21.

\bibitem{BH}
K.~A. Brown and C.~R. Hajarnavis, {\em Homologically homogeneous rings}, Trans.
  Amer. Math. Soc. {\bf 281} (1984), 197--208.

\bibitem{Elkik}
R.~Elkik, {\em Singularit\'es rationnelles et d\'eformations}, Invent. Math.
  {\bf 47} (1978), no.~2, 139--147.

\bibitem{SINGULAR}
G.~M. Greuer, G.~Pfister, and H.~Sch\"onemann, {\em {\sc Singular} 2.0. {A}
  {C}omputer {A}lgebra {S}ystem for {P}olynomial {C}omputations}, Centre for
  Computer Algebra, University of Kaiserslautern, 2001, {\tt
  http://www.singular.uni-kl.de}.

\bibitem{Hochster}
M.~Hochster, {\em Rings of invariants of tori, {C}ohen-{M}acaulay rings
  generated by monomials, and polytopes}, Ann. of Math. (2) {\bf 96} (1972),
  318--337.

\bibitem{KV}
M.~Kapranov and E.~Vasserot, {\em Kleinian singularities, derived categories
  and {H}all algebras}, Math. Ann. {\bf 316} (2000), no.~3, 565--576.

\bibitem{Kawamata}
Y.~Kawamata, {\em Francia's flip and derived categories}, Algebraic geometry
  (Berlin), de Gruyter, Berlin, 2002, pp.~197--215.

\bibitem{Keller1}
B.~Keller, {\em Deriving {DG}-categories}, Ann. Sci. {\'E}cole Norm. Sup. (4)
  {\bf 27} (1994), 63--102.

\bibitem{King}
A.~D. King, {\em Moduli of representations of finite-dimensional algebras},
  Quart. J. Math. Oxford Ser. (2) {\bf 45} (1994), no.~180, 515--530.

\bibitem{Kont1}
M.~Kontsevich, Lecture at Orsay, December 7, 1995.

\bibitem{KuOr}
S.~A. Kuleshov and D.~O. Orlov, {\em Exceptional sheaves on {D}el {P}ezzo
  surfaces}, Izv. Ross. Akad. Nauk Ser. Mat. {\bf 58} (1994), no.~3, 53--87.

\bibitem{Luna}
D.~Luna, {\em Slices \'etales}, Bull. Soc. Math. France {\bf 33} (1973),
  81--105.

\bibitem{Newstead}
P.~E. Newstead, {\em Introduction to moduli problems and orbit spaces},
  Lectures on Mathematics and Physics, vol.~51, Tata Institute of Fundamental
  Research, Bombay, 1978.

\bibitem{reiner}
I.~Reiner, {\em Maximal orders}, London Mathematical Society Monographs. New
  Series, vol.~28, The Clarendon Press Oxford University Press, Oxford, 2003.

\bibitem{RVdB}
I.~Reiten and M.~Van~den Bergh, {\em Tame and maximal orders of finite
  representation type}, Memoirs of the AMS, vol. 408, Amer. Math. Soc., 1989.

\bibitem{Saito}
K.~Saito, {\em Quasihomogene isolierte {S}ingularit\"aten von
  {H}yperfl\"achen}, Invent. Math. {\bf 14} (1971), 123--142.

\bibitem{schofield1}
A.~Schofield, {\em Representations of rings over skew fields}, Lecture Note
  Series, vol.~92, London Mathematical Society, 1985.

\bibitem{VdBSt1}
J.~T. Stafford and M.~Van~den Bergh, {\em Noncommutative resolutions and
  rational singularities}, Michigan Math. J. {\bf 57} (2008), 659--674, Special
  volume in honor of Melvin Hochster.

\bibitem{St1}
R.~P. Stanley, {\em Linear {D}iophantine equations and local cohomology},
  Invent. Math. {\bf 68} (1982), 175--193.

\bibitem{St5}
\bysame, {\em Combinatorics and commutative algebra}, Progress in Mathematics,
  vol.~41, Birkh\"auser Boston Inc., Boston, MA, 1983.

\bibitem{stenstrom}
B.~Stenstr{\"o}m, {\em Rings of quotients}, Die {G}rundlehren der
  mathematischen {W}issenschaften in {E}inzeldarstellungen, vol. 217, Springer
  Verlag, Berlin, 1975.

\bibitem{Steph1}
D.~R. Stephenson, {\em Artin-{S}chelter regular algebras of global dimension
  three}, J. Algebra {\bf 183} (1996), 55--73.

\bibitem{Steph2}
\bysame, {\em Algebras associated to elliptic curves}, Trans. Amer. Math. Soc.
  {\bf 349} (1997), 2317--2340.

\bibitem{VdB1}
M.~Van~den Bergh, {\em {C}ohen-{M}acaulayness of semi-invariants for tori},
  Trans. Amer. Math. Soc. {\bf 336} (1993), no.~2, 557--580.

\bibitem{VdB16}
\bysame, {\em Existence theorems for dualizing complexes over non-commutative
  graded and filtered rings}, J. Algebra (1997), 662--679.

\bibitem{VdB32}
\bysame, {\em Non-commutative crepant resolutions}, The legacy of Niels Henrik
  Abel (Berlin), Springer, Berlin, 2004, pp.~749--770.

\bibitem{VdB31}
\bysame, {\em Three-dimensional flops and noncommutative rings}, Duke Math. J.
  {\bf 122} (2004), no.~3, 423--455.

\bibitem{Watanabe}
K.~Watanabe, {\em Rational singularities with {$k\sp{\ast} $}-action},
  Commutative algebra (Trento, 1981) (New York), Lecture Notes in Pure and
  Appl. Math., vol.~84, Dekker, New York, 1983, pp.~339--351.

\bibitem{YZ1}
A.~Yekutieli and J.~Zhang, {\em Rings with {A}uslander dualizing complexes}, J.
  Algebra {\bf 213} (1999), no.~1, 1--51.

\bibitem{YZ3}
\bysame, {\em Dualizing complexes and perverse sheaves on noncommutative ringed
  schemes}, to appear, 2002.

\end{thebibliography}

\def\cprime{$'$} \def\cprime{$'$} \def\cprime{$'$}
\ifx\undefined\bysame
\newcommand{\bysame}{\leavevmode\hbox to3em{\hrulefill}\,}
\fi

\end{document}